\documentclass[12pt]{article}
\usepackage{amsfonts,amssymb,version}

\begin{document}

\centerline{\Large My encounter with Seshadri}

\centerline{\Large and with the Narasimhan-Seshadri theorem}

\bigskip

\centerline{Nitin Nitsure}

\bigskip

C.S. Seshadri passed away in July 2020. 
The first part of this article contains some reminiscences from the early 1980s, when as a
graduate student in the Tata Institute of Fundamental Research (TIFR), Mumbai,
I had the good fortune to learn from him. 
The second part is on the Narasimhan-Seshadri theorem and my encounter with it.

{\bf I: Some recollections from the early 1980s.}

I joined the School of Mathematics, TIFR, as a research scholar in August 1980.
Seshadri taught us a five or six month long course titled `Geometry-Topology'. It is traditional in
TIFR to have three courses in the first year, titled `Algebra', `Analysis' and `Topology'.
The title `Geometry' was new, and 
so even before the course began, my seniors commented on this departure from the tradition. 

A few words about myself: till joining TIFR, my principal interest was theoretical physics,
but I had taken an MSc degree in Mathematics to escape the physics laboratory
and the hocus-pocus of the lectures.
I had discovered that to learn MSc-level mathematics,
it is enough to sit at home and read by oneself
the various excellent textbooks such as Serge Lang's `Algebra' or Spanier's `Algebraic Topology',
and as a result, I used spend more time in the physics department of the university, and
go to the mathematics department mainly to write the examinations.
When I joined the School of Mathematics of TIFR, my plan was to shift back to physics once I pick up
enough advanced mathematics. I was keen to learn modern geometry, given its
importance in physics, and so I was happy to find that there was going to
be a course on geometry. But I was not so happy with the idea of having to attend lectures -- in 
the previous two years I had not followed any lecture courses,
and had always preferred reading books to attending lectures.
I had not even heard Seshadri's name before.

Seshadri's course was a revelation! He explained later (by which time he had left TIFR)
that the idea was to introduce the basic theory of manifolds, bundles, group actions etc. 
from a Grothendieckian point of view, with emphasis on functorial properties but without
too much abstract machinery, and supplement it with some basic algebraic topology.
In his first lecture, he defined `separated topological spaces' as spaces in which any two
distinct points
have disjoint open neighbourhoods -- which we of course recognized as Hausdorff spaces under a new name
-- and gave various exercises. The first one said that a topological space is separated if and only
if its diagonal $\Delta_X$ is closed in $X\times X$.
The second one said that if $f,g : X\to Y$ are continuous maps where
$Y$ is separated, then the subset defined by $f(x) = g(x)$ is closed in $X$. Yet another one
said that if some point in a topological group is a closed point then the group is separated,
and so on. The whole time was spent is giving about ten `Exercises' followed by about ten more  
`Elementary Exercises'. This baffled us: the elementary exercises looked no easier than the exercises.
That afternoon, as I tried to solve the exercises, I realized that the set
$\{ x\in X\,|\,f(x) = g(x)\}$ is simply the inverse image of the diagonal $\Delta_Y\subset Y\times Y$
under $(f,g) : X\times X \to Y\times Y$, and all the other exercises had similar easy arrow-theoretic
solutions. There was no need to consider points and neighbourhoods, once separatedness was
formulated in terms of the diagonal. I was hooked.   

The second lecture introduced proper maps and proper actions of topological groups.
A map $f: X\to Y$ was defined to be proper if for all $Z$, the induced map
$f\times {\mathop{\rm id}\nolimits}_Z : X\times Z \to Y \times Z$ is closed.
A topological group action $a: X\times G \to X$
was defined to be proper if the induced map $(p_1,a) : X\times G \to X\times X$ is proper. 
Once again, this was followed by a number of
`Elementary Exercises', such as if the action is proper then the quotient $X/G$ is separated. 
I had not seen properness before, and the definition of the properness of an action looked strange
at first sight.
But once again, solving the exercises revealed the meanings of these ideas, and reinforced the power of
arrow-theoretic arguments without recourse to points and neighbourhoods.


\pagestyle{myheadings}
\markright{Nitin Nitsure: A tribute to C S Seshadri.}


Seshadri's lectures were the high point of the week for me for the next few months. 
I said above that Seshadri's course was a revelation. But so was Seshadri's style and persona.
He was not exactly fluent as a speaker, and had to frequently consult his notebook when he wrote
on the board. This gave some much-needed extra time to understand what he was trying to say,
and sometimes (to my great delight) to anticipate what was coming. He was completely relaxed and
clearly in love with the material. The only times when he used to look flustered was
when he arrived to teach the class 30 minutes late instead of his usual 15 minutes late. 
After such extra late arrivals, he would breathlessly say `oh my watch ...', leaving the sentence
incomplete (therefore we could never quite figure out what was the matter with his watch). 
In fact, many mathematical sentences or even whole paragraphs were left incomplete,
with a few suggestive phrases, punctuated by a stream of `you know's! A differential manifold
was defined by drawing two disjoint circles, below them two intersecting circles
and then by performing some kind of a dance with the chalk, showing by tapping on the
board how a point in one of 
the neighbourhoods gets identified with a point in the other neighbourhood.
While performing this dance, he kept looking at us, saying `you know!' and nothing more.
(One of my batch mates was so overcome by this performance 
that he had to leave the room choking, doubled up with laughter.
He ended up as an analytic number theorist -- such is justice.) 
Often, the rest of my day after Seshadri's lecture would go into figuring out what he meant. 
Every lecture was a treat, in which fundamental and beautiful stuff was taken 
out of his magic notebook, and displayed to us. One day, he defined Grassmannians and the tautological
vector bundles on them. The universal property of these bundles
was thrown at us as a challenge, which I could prove only after a hard two or three days of struggle. 
I was anyway used to self-study: what I did not 
know was {\it what} to study. So this course was just made for me! Later when time came for
teaching homology theory, Seshadri (who by then must have developed enough confidence in me)
asked me to lecture in his place, while he sat in the back row. That is when I 
discovered that though he was a woolly speaker, he was a sharp 
listener if and when he chose to pay attention. 

To put all of this in some perspective, the School of Mathematics of TIFR of the early 1980s
was overflowing with brilliant mathematicians of all ages, and had a culture of
unrelenting cleverness. To be regarded as any good, you had to all the time solve tough exercises,
involving imaginative tricks. 
Lectures used to cover huge amounts of material, at breakneck speeds.
The atmosphere was competitive, and it was not uncommon to hear casual evaluatory comments about
anybody not present.
Even in the Canteen conversations, one had to nonstop figure out the puns that were flying about,
or be branded as dim. And here was Seshadri, 
not very fluent, slow to understand jokes (so they said), always very kindly and happy looking, 
without a trace of the biting humor and sarcasm which was the common style, and yet at the very top.
A truly absent minded professor, with goodwill towards all, happy to live in his own
rarefied world of mathematics and Indian classical music. 
This was much before the Chennai Mathematical Institute (CMI) was, as they say, even a gleam in his eye.
Everybody -- including I am sure he himself -- would have been 
astounded if told that he will set up a new world-class institution for mathematical
research and education in Chennai within a few more decades,
struggling against enormous odds, bringing together a team
of dedicated and talented younger people, and enlisting the help of diverse private donors and
government agencies! It must not be just his mathematical eminence, but his entire persona,
that made this possible. 

My frequent interactions with Seshadri as a student ended with his lecture course.
In the final oral examination in April 1981 -- to which I went all prepared to rattle off stuff
such as the proof of the Frobenius theorem on integrable distributions -- he made no mention at all
of the course material, and instead 
asked me a single question: `Have you studied characteristic classes?'.
When I answered that I had not, he asked me to read them, and also read   
Serre's Faiseaux Alg\'ebriques Coh\'erents (FAC) (see [Se]),
as home work in the coming vacation. He told me that it is imperative to acquire French if one is to  
understand advanced mathematics.
That vacation, sitting at home in Pune with a French dictionary in one hand, I struggled
through the early portion of FAC, learning sheaves and abstract varieties.
The basic definitions proposed in FAC were
deeply satisfying to me, something that I have always missed in physics,
which seems to be without clear definitions. 
This made me postpone indefinitely my earlier plans to
return to physics after a year or two in the School of Mathematics, and I ended up spending
almost 40 years there. Not much further interaction with
Seshadri was possible: he was not available for most of the academic year 1981-82 because of a
mathematical emergency
(somebody else's theorem that he had used turned out to have a problem, and he was busy
trying to fix the part that he needed, which he finally accomplished after several months of hard
struggle). In early 1982, he asked me 
to give a series of introductory talks for graduate students
on Hodge theory and Kodaira vanishing to which he came, 
and one day in July 1982, just before he left TIFR for ever, 
he called me to his office and gave me a research problem which was to determine the
singular cohomology of the moduli space of parabolic bundles on a Riemann surface 
(this became my first published mathematics paper). 

I will not narrate the rest of my personal reminiscences of Seshadri here. That is because
the stories that I can tell
-- though they are extraordinary by themselves and are precious to me --
are mostly of the kind that are well known to others.
In the last decade while serving on
the Governing Council of CMI, I made many trips to Chennai, and always looked forward to
meeting him and discussing mathematics, music, history, politics and all other interesting things
under the sun, always coming back energized and more enlightened. 
Many others who were much more closely associated with him have already told their insightful stories.
I will only say that from 1980 to the end of his life, he was
somebody whom I greatly looked up to, and he was always very kind, appreciative, and supportive 
towards me. It is indeed my good fortune to have known such a man and become friends with him,
and his memory will always inspire me.

\bigskip

{\bf II: The Narasimhan-Seshadri Theorem.}

We know that Narasimhan and Seshadri were essentially self-taught as graduate students, who mastered 
huge portions of Seminaire Cartan and other latest such material as it arrived
in TIFR in 1950s. They did this by running a seminar in which they lectured to each other.
Occasional visits to TIFR by eminent mathematicians, and lecture courses by some of these
visitors, gave them inputs about what is important and
what to study. Postdoctoral stints in Paris, mentored by Schwartz (Narasimhan) and Chevalley
(Seshadri) had completed their transition from students to researchers. 
Even before they went to Paris they were aware (thanks to K.G. Ramanathan) of the paper
by Weil [We] in which it is suggested that vector bundles on compact Riemann surface $X$ that are  
associated to unitary representations of the fundamental group $\pi_1(X)$ may have
some unspecified important special properties.
They began their investigation of such bundles after they came back to TIFR ($\sim$1960),
and the result which they found in 1964 is the Narasimhan-Seshadri theorem.

Some of the most important theorems in mathematics can be characterized as `bridge theorems', for they 
connect together two quite different regions. The concept of 
an irreducible unitary representation $\rho : \pi_1(X) \to U(n)$ of the fundamental group of
a compact Riemann surface $X$ arose out of one domain of mathematics, and the concept of a stable
holomorphic vector bundle $X$ arose from a very different domain. The theorem shows that not just
the bundle $E_{\rho}$ on $X$ associated to $\rho$ is a stable
holomorphic vector bundle, but every stable holomorphic vector bundle on $X$ of rank $n$ and degree $0$
so arises up to isomorphism
from exactly one such irreducible unitary representation 
up to conjugacy of representations. There is a similar statement
(though a bit more complicated to state)
for stable holomorphic vector bundles on $X$ of a non-zero degree, which replaces representations
of $\pi_1(X)$ by those representations
of $\pi_1(X - x)$ where $x\in X$ which have a certain scalar local monodromy around $x$. 

The concepts of stability of bundles and stability under group actions 
were defined by Mumford ($\sim$1962) as key points of his Geometric Invariant Theory
(GIT) 
which he used to make a moduli scheme for stable bundles.
All that Narasimhan and Seshadri used from Mumford's work [Mu-2] is the
definition of stability of a bundle, which they knew from
Mumford's talk [Mu-1] in the International Congress for Mathematicians (ICM), Stockholm 1962,
as the volume of the proceedings of the ICM was available in the TIFR library.
The Narasimhan-Seshadri theorem showed what stability means from another perspective, 
and consequently the new field of mathematics -- vector bundles and their moduli spaces 
-- which had been dormant after the early work of Weil, Grothendieck and Atiyah, suddenly became
a rapidly expanding hot luminous star in the mathematical firmament! Over the last 50 years the 
resulting developments have impacted a wide arena, including the Langlands program,
differential geometry, and the new theories of physics -- for example, Witten wrote that     
``The Seshadri-Narasimhan theorem plays a central role in the relationship between conformal
field theory in two dimensions and Chern-Simons gauge theory on three dimensions'' (see [Wi]).

Eminent mathematicians are often feted with birthday conferences. But for the first time
as far as I know, this honor was conferred on a theorem, when the 50th anniversary of the 
Narasimhan-Seshadri Theorem was celebrated with an international conference at CMI in Chennai in 2015.
Both Narasimhan and Seshadri were present and in great form. The Organizing Committee had
assigned to me two jobs: to give an exposition of the Narasimhan-Seshadri theorem as the opening lecture
of the conference (with both Narasimhan and Seshadri present in the audience!),
and to conduct the closing session of felicitatory speeches. 
So probably for the first time after the early 1980s, I tried to read the two papers [N-S-1] and [N-S-2]
of Narasimhan and Seshadri where they prove the theorem.

At this point I should tell the reader that Narasimhan was my thesis supervisor.
After Seshadri left TIFR in 1982, I went on working on my own for some time and 
proved some small results, when Raghunathan pushed me to go and talk with Narasimhan and ask
him to guide me for PhD.
To my greatest good fortune 
Narasimhan agreed, and has 
been my teacher (a permanent, lifelong position in the Indian tradition in its best sense)
from 1984 onwards. This is not the place to write more about Narasimhan, except to say that he  
fully encouraged me to pursue Algebraic Geometry in the Grothendieckian way, even though
his own natural sympathies were much more analytic in those days. 
This is relevant to what follows.

When I looked at those papers in 2015 --
with the advantage of looking at them 50 years later, and with eyes that were trained
by both Narasimhan and Seshadri -- some aspects immediately stood out:\\
(1) Though the theorem is about individual bundles,
the proof is via moduli spaces and maps between them.\\
(2) The moduli space of simple bundles is made in an analytic way using 
the deformation theory of Kodaira and Spencer. \\
(3) The arguments use a whole lot of fundamental 
and beautiful mathematics developed in the 1950s, which is
in the style of Bourbaki but which is pre Grothendieck and pre GIT,
and the formulations are not very functorial. 

Let me elaborate a bit on the above three aspects,
at the cost of repeating things that some readers may
already know.
As Narasimhan explained to me long ago, he and Seshadri were aware of the 
`continuity method' that Klein and Poincar\'e had invented
(independently of each other at about the same time), which
was used to prove the uniformization theorem.
That also is a proof based on moduli-theoretic arguments of a statement that is about individual objects.
Moduli spaces were much in the air in the 1960s, and Seshadri [S-1] had already constructed the Picard
variety in an important case. This is the background of the point (1). 
There is some interesting history related to (2). Namely, during his postdoctoral years, Narasimhan
fell ill and had to spend a few months in a sanatorium in France, where he carried with him
the preprints of the two
fundamental papers of Kodaira and Spencer on deformation theory which were given
to him by Schwartz. He had mastered these while recuperating, and that is how
analytic deformation theory (which was
freshly minted, cutting edge stuff in 1958) became a part of his tool kit. 
In reading the Narasimhan-Seshadri papers, it is exactly the analytic deformation theory part
that caused me difficulty, as I am neither sufficiently well versed in analytic methods, nor
exactly enthusiastic about them, though algebraic deformation theory is something that I love
(see [Ni-5], [Ni-6]).
Now about the point (3): the Narasimhan and Seshadri 
that I first met in the early 1980s were well versed in Grothendieck's philosophy,  
and employed Grothendieck's methods to varying extents. 
If they themselves had to rewrite their proof even 5 years
after they originally wrote it, 
it is reasonable to guess that they would have made it more explicitly functorial and
perhaps even scheme theoretic (for example by using nilpotents). 

So I decided to go all the way in my expository lecture and replace
the Kodaira-Spencer part and the analytic moduli construction part
by the algebraic space of simple bundles made by using the work
of Schlessinger [Sc] ($\sim$1967) on algebraic deformation theory and that of Michael Artin
[Ar] ($\sim$1968) on moduli constructions via algebraic approximation (for pdf slides of the lecture,
see [Ni-7]). 
This was nearly 30 years after my PhD viva examination by Narasimhan and
Seshadri (who was the external examiner), which I remembered when I saw them all ready in the 
first row of the CMI lecture room before my lecture began.
I was quite aware of the enormity of my presumption,  
and felt relieved when all went well --
they liked the reformulation of the proof and they regarded it as both natural and worth the effort.
Both were their usual curious, gracious and generous selves, who wore their own eminence
lightly.

A few years after the Narasimhan-Seshadri theorem, Seshadri introduced the notion of S-equivalence
of semistable bundles (means Jordan-Holder equivalence), and constructed a good moduli space for
semistable bundles, whose points are the S-equivalence classes of semistable
bundles of a given rank and degree (see [S-2]). The term `good moduli space', which has come
into increasing use in the past two decades, indicates certain 
properties of the moduli space of an algebraic stack that originated in 
Seshadri's formulation of the notion of a `good quotient' in the context of GIT.
If one is willing to use the moduli space of Seshadri, then the above proof of the Narasimhan-Seshadri 
theorem further simplifies: the use of Artin's approximation theorem and algebraic spaces 
can be replaced by the use of Seshadri's moduli space. I lectured on this simplified 
proof in the virtual seminar of IIT-B in July 2020
(the video, pdf slides and related notes are available from [Ni-4], [N-5], [Ni-8]).
To my great delight, 
Narasimhan was present in the online audience and seemed to be happy with the presentation
(undoubtedly, he and other experts must have already known that such a proof is possible). 

Seshadri further generalized the ideas from the Narasimhan-Seshadri theorem around 1968 by introducing
the concept of parabolic vector bundles. In collaboration with Vikram Mehta, he constructed a
moduli for these and proved what is known as 
the Mehta-Seshadri theorem [M-S],
which links irreducible unitary representations of the fundamental group of 
a compact Riemann surface with finitely many punctures to stable parabolic bundles
on it. 

Instead of moving in the algebraic direction (which would for instance involve
the Tannakian theory of Deligne and Saavedra),
one can move in the differential geometric
direction, to formulate and prove the Narasimhan-Seshadri theorem in terms of curvature.
In their 1965 paper (see the end of section 10 of [N-S-2]),
Narasimhan and Seshadri cryptically refer to such a formulation over compact Riemann surfaces.
The famous paper [A-B] of Atiyah and Bott in 1980 made the connection between stability
and the Yang-Mills equations.
Donaldson [D-1], [D-2], [D-3] ($\sim$1983-85) gave a new analytic proof of the
Narasimhan-Seshadri theorem  and
also extended the theorem to nonsingular projective varieties of all dimensions
by proving the existence of Hermitian-Einstein metrics on stable bundles.
The work of Hitchin [H] on Higgs bundles ($\sim$1986) is a very important development in this subject,  
that arose from differential geometric ideas, which generalizes the earlier theory
from unitary representations to all linear representations, and vector bundles to Higgs bundles.
It has had unexpected applications, including in the work of Ng\^o on the Fundamental Lemma
(see [Ng] for a survey).
Uhlenbeck, Yau, Corlette and Simpson (see [U-Y-1], [U-Y-2], [C], [Si-1], [Si-2], [Si-3]) 
made further advances in the differential geometric direction inspired by the
Narasimhan-Seshadri theorem. I was also able to make a modest contribution to the
developing theory (see [Ni-3]). Many distinguished mathematicians
have told me that one of the first papers that they were asked to read when they began their graduate
study (at some of the top places in the world) was the 1965 Narasimhan-Seshadri paper, and it
played an important role in shaping their outlook.
Today, the Narasimhan-Seshadri theorem is tightly woven
into the very fabric of a vast domain of algebraic geometry, 
and it is a fair guess that this is going
to remain so for a very very long time to come.

The last time I met Seshadri was when I made a condolence visit
to his house in November 2019 after the death of Sundari, his wife.
He talked about Sundari -- I had never heard him speak about her before, and
saw what a pillar of strength she must have been for him. Then changing the subject, he asked
me to explain my recent theorem [Ni-9] on quadrilateral gaps for affine connections on manifolds.
We discussed mathematics and many other things, and laughed like we always did,
but with a heavy heart. When he came to the gate to see me off, 
both of us were aware of the passing of time. It does not stop for anybody.

My salute!

\newpage

\parskip 3pt

{\bf References}

[Ar] Artin, M. : Algebraization of formal moduli I. In {\it Global analysis,
Papers in Honor of K. Kodaira}, 21-71. Tokyo Univ press 1969. 

[A-B] Atiyah, M.F. and Bott, R. : Yang-Mills on Riemann Surfaces.
Phil. Trans. R. Soc. Lond. 308 (1983) 523-615.

[C] Corlette, K. : Flat $G$-bundles with canonical metrics. J. Differential Geom. 28 (1988) 361-382.

[D-1] Donaldson, S.K. : A new proof of a theorem of Narasimhan and Seshadri.
J. Differential Geom. 18 (1983) 269-277. 

[D-2] Donaldson, S. K. :
Anti self-dual Yang-Mills connections over complex algebraic surfaces and stable vector bundles.
Proc. London Math. Soc. (3) 50 (1985) 1-26.

[D-3] Donaldson, S. K. : Infinite determinants, stable bundles and curvature.
Duke Math. J. 54 (1987), 231–247. 

[H]  Hitchin, N. J. : The self-duality equations on a Riemann surface.
Proc. London Math. Soc. (3) 55 (1987) 59–126. 

[M-S] Mehta, V.B. and Seshadri, C.S. :
Moduli of vector bundles on curves with parabolic structures. Math. Ann. 248 (1980) 205-239.

[Mu-1] Mumford, D. : Projective invariants of projective structures and applications.
Proc. Internat. Congr. Mathematicians (Stockholm, 1962) 526-530, Inst. Mittag-Leffler 1963. 

[Mu-2] Mumford, D. : Geometric invariant theory. Springer 1965. 

[N-S-1] Narasimhan, M.S. and Seshadri, C.S. : Holomorphic vector bundles on a compact Riemann surface.
Math. Annalen 155 (1964) 69-80.

[N-S-2] Narasimhan, M.S. and Seshadri, C.S. :
Stable and unitary vector bundles on a compact Riemann surface.
Annals of Math. 82 (1965) 540-567.

[Ng] Ng\^o  B.C.  : Survey on the fundamental lemma. Current developments in mathematics 2009,
1-21, Int. Press 2010. \\
https://www.math.uchicago.edu/~ngo/survey.pdf

[Ni-1] Nitsure, N. : Cohomology of the moduli of parabolic vector bundles.
Proc. Indian Acad. Sci. Math. Sci. 95 (1986) 61-77.

[Ni-2] Nitsure, N. : Parabolic bundles. A tribute to C. S. Seshadri (Chennai, 2002), 28-33,
Trends Math., Birkhäuser 2003.  

[Ni-3] Nitsure, N. : Moduli of semistable pairs on a curve.
Proc. London Math. Soc. (3) {\bf 62} (1991) 275-300.

[Ni-4] Nitsure, N. : Vector Bundles on Compact Riemann Surfaces. Lecture notes, 
NCM summer school IIT-M 2007. \\
https://www.ncmath.org/lecture-notes

[Ni-5] Nitsure, N. : Deformation Theory for Vector Bundles. Updated version (2011) of
an earlier article in vol 359 LMS lecture notes in honour of Peter Newstead. \\
https://www.newton.ac.uk/files/seminar/20110105153016301-152666.pdf

[Ni-6] Nitsure, N. : Deformation theory and Moduli Spaces. Masterclass 2014,
Center for Quantum Geometry of Moduli, Aarhus. \\
Slides: https://qgm.au.dk/fileadmin/www.qgm.au.dk/Events/2014/Slides$\_$Nitin.pdf\\
Videos: https://qgm.au.dk/video/mc/triple-masterclass/index.html

[Ni-7] Nitsure, N. : Slides of the lecture in the 50th anniversary conference of
the Narasimhan-Seshadri Theorem at CMI (2015). \\
https://www.cmi.ac.in/$\sim$ramadas/NS@50/Nitsure.pdf

[Ni-8] Nitsure, N. : Slides and videos of two lectures on 
the Narasimhan-Seshadri Theorem in the online mathematics seminar of IIT-B (2020). \\
http://www.math.iitb.ac.in/$\sim$ronnie/past-talks.html

[Ni-9]  Nitsure, N. : Curvature, torsion and the quadrilateral gaps. arXiv:1910.06615
To appear in Proc. Indian Acad. Sci. Math. Sci. 131 (2021). 

[Se] Serre, J.-P. : Faiseaux Alg\'ebriques Coh\'erents. Annals of Math. 61 (1955) 197-278.

[S-1] Seshadri, C.S. : Vari\'et\'e de Picard d'une vari\'et\'e compl\`ete. Ann. Mat. Pura Appl.
(4) 57 (1962), 117-142.

[S-2] Seshadri, C.S. : Space of unitary vector bundles on a compact Riemann surface.
Ann. of Math. (2) 85 (1967) 303-336.

[Sc] Schlessinger, M. : Functors of Artin rings. Trans. A.M.S. 130 (1968) 208-222.

[Si-1] Simpson, C.T. : Harmonic bundles on noncompact curves.
J. Amer. Math. Soc. 3 (1990), no. 3, 713-770. 

[Si-2] Simpson, C.T. : Higgs bundles and local systems. IH\'ES Sci. Publ. Math. 75 (1992).

[Si-3] Simpson, C.T. : Moduli of representations of the fundamental group of a smooth
projective variety - I and II. IH\'ES Publ. Math. 79 (1994) 47-129 and 
80 (1994) 5-79. 

[U-Y-1] Uhlenbeck, K. and Yau, S.T. : On the existence of Hermitian-Yang-Mills connections
in stable vector bundles. Comm. Pure Appl. Math. 39 (1986), 257-293.

[U-Y-2] Uhlenbeck, K. and Yau, S.T. : 
A note on our previous paper: "On the existence of Hermitian-Yang-Mills connections
in stable vector bundles''. Comm. Pure Appl. Math. 42 (1989) 703-707. 

[We] Weil, A. : G\'en\'eralisation des fonctions ab\'eliennes. J. Math. pures et appl. 17 (1938) 47-87. 

[Wi] Witten, E : @witten271 reply to @raghumahajan and @scroll$\_$in, Twitter 2020. \\
https://twitter.com/raghumahajan/status/1289453265868828673?lang=en

\bigskip

\footnotesize

Nitin Nitsure \\
Retired Professor, School of Mathematics \\
Tata Institute of Fundamental Research \\
Mumbai 400 005, India \\
email: nitsure@gmail.com

\end{document}